\documentclass{article}
 \usepackage{yhmath}
 \usepackage{amssymb}
 \usepackage{float}
\usepackage{amsmath} 

\newcommand{\R}{\mathbb{R}}
\newcommand{\N}{\mathbb{N}}
\newtheorem{defn}{Definition}

\newtheorem{rem}{Remark}
\newtheorem{cor}{Corollary}
\newtheorem{lem}{Lemma}
\newtheorem{theo}{Theorem}

\newtheorem{Ex}{Example}

\usepackage[english]{babel}

\usepackage[letterpaper,top=2cm,bottom=2cm,left=3cm,right=3cm,marginparwidth=1.75cm]{geometry}

\usepackage{amsmath}
\usepackage{graphicx}
\usepackage[colorlinks=true, allcolors=blue]{hyperref}
\title{On the Failure of the Łojasiewicz Property for Smooth Ideals with Dense Regular Loci}
\author{Abdelhafed El Khadiri}
\begin{document}

\maketitle
\protect\footnote{Math Subject Classification (MSC 2020) Primary: 26B05, 58A20, Secondary: 58C25, 57N80.}

\begin{abstract}
Let $\Omega \subset \mathbb{R}^n$ be an open set, and let $\mathcal{E}(\Omega)$ denote the ring of smooth functions on $\Omega$. For a finitely generated ideal $I \subset \mathcal{E}(\Omega)$, we denote by $Z(I)$ its zero set. A classical result of René Thom asserts that if $I$ is a Łojasiewicz ideal, then $Z(I)$ contains an open dense subset of regular points.

In this paper, we investigate the converse problem from a geometric perspective: does the existence of an open dense regular locus in $Z(I)$ imply that $I$ is a Łojasiewicz ideal? We show that this implication fails in general, and we identify a genuinely singular obstruction responsible for this failure.

More precisely, the obstruction arises from the accumulation of infinitely many smooth branches whose tangential behavior degenerates at a singular point. This phenomenon is illustrated by a construction inspired by the classical Hawaiian earring, namely the union of infinitely many pairwise tangent circles accumulating at a common point. We prove that any smooth function whose zero set coincides with this space must be flat at the accumulation point.

These results highlight a sharp distinction between analytic and smooth singularities, and demonstrate that the topology of the regular locus alone is insufficient to ensure Łojasiewicz-type properties in the $C^\infty$ category.
\end{abstract}
\section{Introduction}

The study of zero sets of smooth functions and ideals in the ring 
$\mathcal{E}(\Omega)$ of infinitely differentiable  functions on an open set 
$\Omega \subset \mathbb{R}^n$ lies at the intersection of real analytic geometry, 
singularity theory, and differential topology. 
A central tool in this context is the {\L}ojasiewicz inequality, which provides a 
quantitative relation between a function and its gradient, and plays a fundamental 
role in the analysis of subanalytic sets, stratifications, and resolution of singularities.

A finitely generated ideal $I = (f_1, \ldots, f_k)\mathcal{E}(\Omega)$ is said
to be a {\L}ojasiewicz ideal if there exists a function $g \in I$ that satisfies
a {\L}ojasiewicz inequality with respect to
$
Z(I) := \{ x \in \Omega \mid f_1(x) = f_2(x) = \cdots = f_k(x) = 0 \}.
$
See the definition below.\\
A classical result of Ren\'e Thom \cite{thom1967some} asserts that if $I$ is a finitely generated  Lojasiewicz ideal, 
then its zero set $Z(I)$ contains an open dense subset of smooth points.\\  In analytic geometry this is expected: analytic ideals naturally enforce 
geometric regularity on their zero sets.  
However, in the smooth category such regularity properties are more subtle, since 
arbitrary smooth functions may exhibit behavior impossible in the analytic setting.

This raises a natural question: 
\emph{does the converse hold?}  
Does the presence of an open dense set of smooth points in $Z(I)$ imply that 
the ideal $I$ is \L{}ojasiewicz?  
The purpose of the present work is to examine this question and to identify 
mechanisms by which such a converse can fail.\\
To illustrate this phenomenon, we focus on a classical topological example:  
the \emph{Hawaiian earring}, defined as the union of the circles
\[
C_n = \left\{ (x_1,x_2) \in \mathbb{R}^2 : 
\left(x_1 - \frac{1}{n}\right)^2 + x_2^2 = \frac{1}{n^2} \right\}, \qquad n \ge 1.
\]
This compact set consists of infinitely many circles tangent at the origin and 
accumulating there.  
It is well known that the Hawaiian earring has highly pathological local 
topology at the origin: it is not locally contractible, not semianalytic, and 
cannot arise as the zero set of any nontrivial real analytic function.

We consider smooth functions with real values $f \in \mathcal{E}^\infty(\R^2)$ whose zero 
set $Z(f)$ is \emph{exactly} the Hawaiian earring.  
Such functions exist by Whitney's extension theorem, see theorem 3.1, Ch IV, \cite{Tougeron1972}, yet they exhibit extremely 
degenerate behavior at the origin.  
Our first main result states that any such function must vanish to infinite order.
\begin{theo}
Let $f \in \mathcal{E}^\infty(\R^2)$ be such that $Z(f)$ is the Hawaiian earring. Then $f$ is flat at the origin:
\[
\frac{\partial^{\alpha +\beta}f}{\partial x_1^\alpha \partial x_2^\beta}(0,0)
 = 0 
\qquad \text{for all } \alpha,\beta \in\N.
\]
\end{theo}
The proof relies on purely geometric considerations: any nonzero term of the 
Taylor expansion of $f$ at $(0,0)$ would define an algebraic curve of finite 
order, and no such curve can contain infinitely many smooth arcs with tangency 
and curvature behavior comparable to the circles $C_n$.  
Thus the geometry of the Hawaiian earring forces infinite-order vanishing.

This result demonstrates that certain pathological smooth zero sets necessarily 
force infinite-order degeneracy, and therefore cannot arise from \L{}ojasiewicz 
ideals.  
In particular, the ideal generated by a function defining the Hawaiian earring is 
never \L{}ojasiewicz, despite the fact that its zero set is smooth away from a 
single point.

\section{Preliminaries on \L{}ojasiewicz Ideals}\label{sec:preliminaries}
In this section we recall the basic definitions and properties of 
\L{}ojasiewicz inequalities and \L{}ojasiewicz ideals in the setting of smooth 
functions.  Our presentation follows the classical framework developed by 
\L{}ojasiewicz \cite{Lojasiewicz1965}, Malgrange \cite{Malgrange1966},   Tougeron \cite{Tougeron1972}, and subsequent refinements 
in real-analytic and subanalytic geometry.

\subsection{\L{}ojasiewicz ideals}
Let $\Omega$ be an open subset of $\mathbb{R}^n$, and let $\mathcal{E}^\infty(\Omega)$ 
denote the Fréchet algebra of smooth functions on $\Omega$.
Let $X$ be a closed subset of $\Omega$. 
An element $\varphi \in \mathcal{E}^\infty(\Omega)$ is said to satisfy 
the {\L}ojasiewicz inequality with respect to $X$ if, for every compact subset 
$K \subset \Omega$, there exist constants $C > 0$ and $\nu \ge 0$ such that, 
for every $x \in K$, we have
\[
|\varphi(x)| \ge C\, d(x,X)^\nu .
\]
\begin{defn}
A finitely generated ideal $I = (f_1, \ldots, f_k)\mathcal{E}^\infty(\Omega)$ 
is said to be a {\L}ojasiewicz ideal if there exists $g \in I$ that satisfies 
the {\L}ojasiewicz inequality with respect to 
\[
Z(I) := \{ x \in \Omega \mid f_1(x)=\cdots=f_k(x)=0 \}.
\]
\end{defn}

\begin{rem}
In this case, for any system of generators $g_1,\ldots,g_p$ of $I$, 
the functions $\sum\limits_{j=1}^p g_j^2$ and $\sum\limits_{j=1}^p |g_j|$ both satisfy 
the {\L}ojasiewicz inequality with respect to $Z(I)$.
\end{rem}
It should be noted that the property of an ideal to be {\L}ojasiewicz is a 
local one: if an ideal is {\L}ojasiewicz on an open set $U$, then the induced 
ideal on any smaller open subset is also {\L}ojasiewicz.\\
Any analytic ideal is {\L}ojasiewicz, as this follows from the fundamental 
{\L}ojasiewicz inequality for analytic functions \cite{Lojasiewicz1965}. Moreover, any finitely 
generated ideal that is closed in $\mathcal{E}^\infty(\Omega)$ is also a 
{\L}ojasiewicz ideal, see Corollary 4.4, Ch V, \cite{Tougeron1972}. Indeed, by Whitney's spectral theorem, a smooth function 
belongs to a closed finitely generated ideal $I=(f_1,\ldots,f_k)$ if and only if, 
for every point $x\in\Omega$, its Taylor expansion at $x$ lies in the ideal 
generated by the Taylor expansions $T_x f_1,\ldots,T_x f_k$ in the formal power 
series algebra $\mathbb{R}[[X-x]]$. 
However, it should be emphasized that there exist {\L}ojasiewicz ideals 
which are not closed, in dimension $n >1$, see Examples 4.8, Ch V, \cite{Tougeron1972}.
In dimension one, closed finitely generated ideals are exactly the finitely generated Łojasiewicz ideals. We include a proof of this equivalence in  section 7.
\subsection{Smooth points of  the locus of zeros}
\begin{defn}[Smooth point]
Let $E \subset \mathbb{R}^n$ be closed set. 
A point $x \in E$ is said to be smooth if, in a neighborhood of $x$, 
the set $E$ coincides with a $k$-dimensional embedded submanifold of $\mathbb{R}^n$, 
where $k$ is its local dimension at $x$ and may vary with $x$.
\end{defn}
\begin{Ex}
 The \emph{Hawaiian earring}, 
 defined as the union of circles 
\[
\bigcup_{n=1}^{\infty} \left\{(x,y)\in \mathbb{R}^2 : \left(x-\frac{1}{n}\right)^2 + y^2 = \frac{1}{n^2} \right\},
\] 
is smooth at every point except the origin, where infinitely many circles accumulate.  
Hence, the set of smooth points of the Hawaiian earring is dense in $E= \bigcup_{n=1}^{\infty} \left\{(x,y)\in \mathbb{R}^2 : \left(x-\frac{1}{n}\right)^2 + y^2 = \frac{1}{n^2} \right\}$.
\end{Ex} 
\begin{figure}[!h]
    \centering
\includegraphics[width=0.5\linewidth]{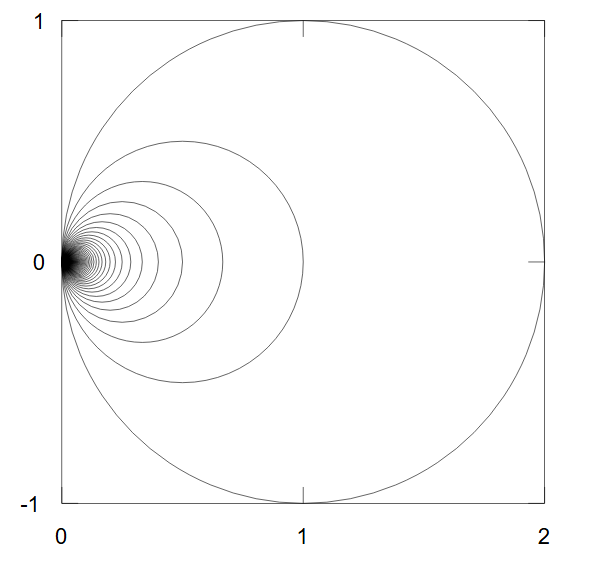}
\caption{$Hawaiian_ Earring$}
\label{fig:placeholder}
\end{figure}
Theorem of Rene Thom states:
\begin{theo}[\cite{thom1967some}]
If $I \subset \mathcal{E}^\infty(\Omega)$ is a \L{}ojasiewicz ideal, then 
$Z(I)$ contains an open dense subset of smooth points.
\end{theo}

Thus \L{}ojasiewicz ideals impose strong geometric regularity on their zero 
sets.  A central theme of the present paper is that the \emph{converse 
statement fails}, even in the simplest settings of functions in two variables.

These preliminary facts form the backdrop for the geometric analysis undertaken 
in later sections, where we show that the Hawaiian earring forces infinite-order 
flatness for any defining function, and we derive a degenerate 
\L{}ojasiewicz inequality adapted to such configurations.
\section{Flatness of indefinitely differentiable functions vanishing on the Hawaiian earring}

\begin{theo}
Let $f:\mathbb{R}^2\to\mathbb{R}$ be a $C^\infty$ function.  
Assume that $f$ vanishes on the Hawaiian earring
\[
\mathcal H=\bigcup_{n\ge1} C_n, \qquad 
C_n=\left\{(x,y):\left(x-\tfrac1n\right)^2+y^2=\tfrac1{n^2}\right\}.
\]
Then $f$ is flat at the origin:    
\[
\frac{\partial^{\alpha_1 +\alpha_2 f}}{\partial x_1^{\alpha_1} \partial x^{\alpha_2}_2}(0,0)
 = 0 
\qquad \text{for all } \alpha =(\alpha_1,\alpha_2) \in\N^2.
\]
\end{theo}

{\bf{Proof.}}
If $f$ is identically zero in a neighborhood of $(0,0)$ there is nothing to prove.
Otherwise, let $p\ge0$ be the smallest integer such that the Taylor expansion of 
$f$ at the origin has a nonzero homogeneous component of degree $p$.  
Thus we may write, for $(x,y)$ near $(0,0)$,
\[
f(x,y)=P_p(x,y)+R_p(x,y),
\]
where $P_p$ is a nonzero homogeneous polynomial of degree $p$, and
$$R_p(x,y)=o(\|(x,y)\|^p)\,\,\,\,\,\mbox{as} \,\,\,\,\, (x,y)\to(0,0).$$
Let $u=(u_1,u_2)$ be a unit vector with $u_1\neq 0$ (a nonvertical direction).  
For each $n\ge1$, the ray $\{tu : t\ge 0\}$ intersects the circle $C_n$ 
at a point of the form $t_n(u)u$ with $t_n(u)>0$, and $t_n(u)\to 0$ as $n\to\infty$.  
Since $f$ vanishes on $C_n$, we have
\[
0 = f(t_n(u)u)=t_n(u)^p P_p(u) + R_p(t_n(u)u).
\]
Dividing by $t_n(u)^p$ yields
\[
0 = P_p(u) + \frac{R_p(t_n(u)u)}{t_n(u)^p}.
\]
Because $R_p=o(\|(x,y)\|^p)$, the second term tends to $0$ as $n\to\infty$.  
Thus
\[
P_p(u)=0.
\]

Hence $P_p$ vanishes on every unit vector with $u_1\neq 0$.  
This set is an open arc of the unit circle, so by continuity $P_p$ vanishes on a nontrivial arc.  
A nonzero homogeneous polynomial cannot vanish on an open arc of the unit circle.  
Therefore $P_p\equiv 0$, contradicting the minimality of $p$.

Thus there is no nonzero homogeneous term in the Taylor expansion of $f$ at $(0,0)$,  
so all derivatives of $f$ at the origin vanish.  
Hence $f$ is flat at $(0,0)$.
\section{Łojasiewicz Inequality and the Hawaiian Earring}
\begin{theo}
Let $f \in C^\infty(\mathbb{R}^2)$ and assume its zero set is the Hawaiian earring
\[
\mathcal H:= \bigcup_{k=1}^{\infty} \left\{ (x,y) \in \mathbb{R}^2 : \left(x - \frac{1}{k}\right)^2 + y^2 = \frac{1}{k^2} \right\}.
\]
Then $f$ is flat at the origin and does \emph{not} satisfy a Łojasiewicz inequality with respect to $H$ at $0$; that is, there exist no constants $C>0$ and $\theta>0$ such that
\[
|f(x)| \ge C \, \mathrm{dist}(x,H)^\theta
\]
for all $x$ in a neighborhood of $0$.
\end{theo}
{\bf{Proof.}}
Since $f \in C^\infty$ vanishes on the Hawaiian earring, which has tangent directions at $0$ dense in the unit circle $S^1$, Theorem 4 implies that $f$ is flat at the origin:
\[
D^\alpha f(0) = 0 \quad \text{for all multi-indices } \alpha.
\]
Let $x \in \mathbb{R}^2$ be close to $0$. Consider the sequence of circles in $H$ with radii $r_k = 1/k$ and centers $c_k = (1/k,0)$.  
Choosing $k \sim 1/\|x\|$, the circle $C_k$ closest to $x$ has radius $r_k \sim \|x\|$.  
Hence, for points near $0$,
\[
\mathrm{dist}(x,H) \sim \|x\|.
\]
Assume, for contradiction, that there exist constants $C>0$ and $\theta>0$ such that
\[
|f(x)| \ge C \, \mathrm{dist}(x,H)^\theta
\]
for all $x$ near $0$.  

Since $f$ is flat at $0$, for any integer $N>0$, there exists a neighborhood $U$ of $0$ such that
\[
|f(x)| \le \|x\|^N \quad \text{for all } x \in U.
\]

But $\mathrm{dist}(x,H) \sim \|x\|$, so
\[
|f(x)| \le \|x\|^N \sim \mathrm{dist}(x,H)^N.
\]

Choosing $N > \theta$ gives
\[
|f(x)| \ll \mathrm{dist}(x,H)^\theta
\]
for $x$ sufficiently close to $0$, which contradicts the assumed Łojasiewicz inequality.  
Therefore, the ideal generated by
 the function $f$ in the the ring  $\mathcal{E}^\infty(\R^2)$ is not a 
 Lojasiewicz ideal.

\section{ Finitely generated Lojasiewicz ideal   in dimension one.}
\begin{defn}
    Let $\Omega\subset \R^1$ be an open set. Consider $f\in\mathcal{E}^\infty(\Omega)$ and put $Z(f)=\{x\in\Omega\,/\, f(x)= 0\}$. We say that $x\in Z(f)$ is a zero  of finite order, if there exists $k\in\N^*$ such that $ f^{k}(x)\neq 0$. (The function $f$ is not flat at $a$).
\end{defn}
\begin{theo}
 Let $f\in\mathcal{E}^\infty(\Omega)\setminus \{0\}$, such   that  $Z(f)\neq \emptyset$ and satisfies de  $L$ojasiewicz 
 inequality. Then all the zeros of the function $f $ are of finite order.
\end{theo}

Before we begin the proof of the theorem, we will need an elementary result about sequences. We prefer to state and prove it first.
\begin{lem}
    Let $(x_n)_n$ be a decreasing sequence of real numbers such that $\lim\limits_{n\to \infty} x_n = 0$. We define the sequence $ v_n = x_n - x_{n+1} \geq 0$. Then there exist a decreasing subsequence of $(v_n)_n$.
\end{lem}
{\bf{Proof of the lemma.}}
We construct a strictly decreasing subsequence \((v_{n_k})\) as follows:\\
Choose \(n_1\) such that \(v_{n_1} < 1\).
 Suppose we have already chosen indices \(n_1 < n_2 < \dots < n_k\) such that:
    \[
    v_{n_1} \ge v_{n_2} \ge \dots \ge v_{n_k}.
    \]
    Since \(v_n \to 0\), there exists an index \(n_{k+1} > n_k\) such that:
    \[
    v_{n_{k+1}} < \min\left(v_{n_k}, \frac{1}{k+1}\right).
    \]
 This guarantees that:
    \[
    v_{n_1} > v_{n_2} > \dots > v_{n_k} > v_{n_{k+1}} > \dots
    \]
Thus, we have constructed a strictly decreasing subsequence \((v_{n_k})\) of the original. 

{\bf{Proof of the theorem.}}
We suppose that $f$ has a zero $a\in Z(f)$ which is not of finite order. Then,  for all $m\in \N,\,\,\,\, f^{(m)}(a) =0$. Since $f$ is not identically nulle, there exists $ b\in \Omega$ with $f(b)\neq 0$. We can suppose that $ a < b$ and $a$ is the closest to $b$. Let $\varepsilon >0$, such that $ a+ \varepsilon < b$, so $f$ admits a zero in 
$]a, a+\varepsilon[$. Indeed, suppose $f$ has no zeros in $]a, a+\varepsilon[$. By Taylor's theorem, for each $k\in \N$ there exists $C_k >0$ such that:
$$ \forall x\in ]a, a+\varepsilon[,\,\,\,\,\,|f(x)|\leq  C_k|x-a|^k = C_k d(x, Z(f))^k,$$ which contradicts the $L$ojasiewicz inequality. Hence $a$ is a limit point of zeros of $f$ in $]a, b[  $. The zeros of $f$
in $]a, b[$ can have no other limit points, for such a limit point would be a zero of infinite order of $f$, contradicting the choice of $a$.
Therefore the zeros of $f$  in $]a, b[$  are a
decreasing sequence $(z_k)_{k\in\N}$  with $z_k\to a$.\\
By Lemma 9,  if we replace the sequence $(z_k)_{k\in\N}$ by a subsequence, we can assume  that the sequence $ S_k = (z_k - z_{k+1})$ is decreasing.\\ 
Since the function $f$ satisfies the $L$ojasiewicz inequality. For the compact set $[a,b]$ there exist $C>0$ and $\alpha \geq 0$ such ,
\begin{equation}
 |f(x)|\geq C d(x, Z(f))^\alpha, \,\,\,\,\,\,\mbox{for all } \,\,\,\,\,\,x\in[a,b].   
\end{equation}
let q be the smallest integer greater than or equal to $\alpha$  and put
$C_1 = \sup\limits_{x\in[a,b]}d(x, Z(f))$. We have 
$$|f(x)|\geq C C_1^\alpha \left(\frac{d(x, Z(f))}{C_1}\right)^\alpha \geq \frac{C}{ C_1^{q-\alpha}}d(x, Z(f))^q, \,\,\,\,\,\,\mbox{for all } \,\,\,\,\,\,x\in[a,b].$$
We then see that we can assume in inequality (2.2) that $\alpha$ is an integer greater than or equal to $0$.\\
Note that the integer $\alpha$  in (2.2) is different from $0 $ because $f(a)=0$, hence $\alpha \geq 1$.

Suppose $\alpha =1$ and put $ y_k= \frac{z_k+z_{k+1}}{2}$. Expanding $f$
about $z_k$ by Taylor's theorem up to order zero, we obtain
$$ C d(y_k,Z(f)) = C|y_k - z_k|\leq |f(y_k)| = |y_k - z_k| |f'(u_k)|,$$
where $u_k\in ]y_k, z_k[$. Hence, we have $ 0 < C \leq |f'(u_k)|$.
Letting $k\to \infty$, we obtain $ 0 < C \leq |f'(a)|$, contradicting the fact that $f$ has a zero of infinite order at $a$.

For the sake of clarity, let us suppose that $\alpha =2$, before dealing with the general case.
 
 Expanding $f$
about $z_k$ by Taylor's theorem up to order one, we obtain
\begin{equation}
   C d(y_k,Z(f))^2 = C|y_k - z_k|^2\leq |f(y_k)| = |y_k - z_k| |f'(z_k)| + \frac{|y_k - z_k|^2}{2}|f''(u_k)|, 
\end{equation}
where $u_k\in ]y_k, z_k[$. \\ 
Since $f(z_{k+1}) = f(z_k) = 0$, there exists $x^1_k\in] z_{k+1}, z_k[$ such that $f'(x^1_k) =0$.\\
Expanding $f'$ about $x^1_k$, we obtain
$$ f'(z_k) = f'(x^1_k) + (z_k -x^1_k) f^{"}(r_k),$$
where $x^1_k < r_k < z_k $. Since $f'(x^1_k) = 0$ and $|z_k - x^1_k|\leq 2|y_k - z_k |$, we have
$$|f'(z_k)|\leq 2|y_k -z_k| .|f''(r_k)|.$$
By (2.3), we obtain
$$ C|y_k - z_k|^2 \leq 2 |y_k - z_k|^2|f''(r_k)| + \frac{|y_k - z_k|^2}{2}|f''(u_k)|.$$
Therefore, we have
$$ 0 <C \leq 2 |f''(r_k)| + \frac{1}{2}|f''(u_k)|,$$
where $u_k\in ]y_k, z_k[$ and $x^1_k < r_k < z_k $. We remark that $ x^1_k, u_k,r_k,\in]z_{k+1}, z_k[$.\\
Letting $k\to \infty$, we obtain 
$$ 0 <C \leq 2 |f''(a)| + \frac{1}{2}|f''(a)|,$$
contradicting the fact that $f$ has a zero of infinite order at $a$.

Suppose that in (2.2), $\alpha = m > 1$. Expanding  $f$ about $z_k$ by Taylor's theorem up to order $m-1$, we obtain
\begin{align*}
C d(y_k,Z(f))^m = C|y_k - z_k|^m &\leq |f(y_k)|\\
 &\leq \sum\limits_{j=1}^{m-1}\frac{1}{j!} |y_k - z_k|^j |f^{(j)}(z_k)| + \frac{|y_k - z_k|^m}{m!}|f^{(m)}(v_m)|, 
\end{align*}
where $ v_m\in ]y_k,z_k[$.\\
We can choose
\[
  \begin{array}{ccc}
    x_k^1\in ]z_{k+1},z_k[\,\,\,\,\,\,\mbox{with}\,\,\,\,\,\, f'(x_k^1) = 0
    \\
    x_k^2\in ]z_{k+2},z_k[\,\,\,\,\,\,\mbox{with}\,\,\,\,\,\, f''(x_k^2) = 0
    \\
    x_k^3\in ]z_{k+3},z_k[\,\,\,\,\,\,\mbox{with}\,\,\,\,\,\, f'''(x_k^3) = 0\\ 
    \vdots & &  \\
    x_k^{m-1}\in ]z_{k+m-1},z_k[\,\,\,\,\,\,\mbox{with}\,\,\,\,\,\, f^{(m-1)}(x_k^{m-1}) = 0
  \end{array}
  \]
  Expanding $f'$ about $x^1_k$, we obtain, for each $k\in \N$,
$$ f'(z_k) = f'(x^1_k) + (z_k -x^1_k) f^{"}(r_k),$$
where $x^1_k < r_k < z_k $. Since $f'(x^1_k) = 0$ and $|z_k - x^1_k|\leq 2|y_k - z_k |$, we have
$$|f'(z_k)|\leq 2|y_k -z_k| .|f''(r_k)|.$$
Expanding $f''$ about $x^2_k$, we obtain, for each $k\in \N$,
$$ |f''(r_k)| =   |z_k -x^2_k||f^{"}(s_k)|,$$
where $x^2_k < s_k < z_k $. We have 
$$ |z_k - x_k^2|\leq |z_k - z_{k+2}|\leq  |z_k - z_{k+1}| + |z_{k+1} - z_{k+2}|.$$
Since the sequence $ S_k = (z_k - z_{k+1})$ is decreasing, we have
$$ |z_k - x_k^2|\leq 2 |z_k - z_{k+1}| = 4 |y_k - z_k|.$$
Continuing in this manner (the last step consists of expanding $f^{(m-1)}$ about $x_k^{m-1}$), we find 
$$ C|z_k - y_k|^m\leq C_1|z_k - y_k|^m |f^{(m)}(u_k)| +\ldots + C_m|z_k - y_k|^m |f^{(m)}(v_k)|,$$
where $z_{k+m-1} < u_k,\ldots\ldots, v_k < z_k$. Therefore, for all $k\in\N$,
$$ 0 < C \leq C_1 |f^{(m)}(u_k)| +\ldots\ldots + C_m |f^{(m)}(v_k)|.$$
Clearly $u_k\to a,\ldots\ldots, v_k\to a,\,\,\,\,\,\mbox{as}\,\,\,\,\, k\to \infty$.  Therefore
$$ 0 < C \leq C_1 |f^{(m)}(a)| +\ldots\ldots +C_m |f^{(m)}(a)|,$$
contradicting the fact that $f$ has a zero of infinite order at $a$. Since $f$ is not identically zero, the function $f$ has only zeros of finite order, which proves the theorem.
\begin{cor}
 Let $f\in\mathcal{E}^\infty(\Omega)\setminus \{0\}$, ( $\Omega\subset \R$ is an open set ), such   that  $Z(f)\neq \emptyset$ and satisfies de  $L$ojasiewicz 
 inequality.  Then each $a\in Z(f)$ is an isolated point in $Z(f)$, that is,
 there exists an open set $U\subset\Omega$ such that $U\cap Z(f) = \{a\}$.
\end{cor}

We will prove that any finitely generated ideal in the ring of smooth real-valued functions on an open set $\Omega \subset \R$,   \( \mathcal{E}^\infty(\Omega) \), whose common zero set is discrete, is principal. The proof uses the fact that locally such ideals are principal and that a smooth partition of unity allows one to patch these local generators into a global one.
\begin{theo}
    Let \( I \subset \mathcal{E}^\infty(\Omega) \) be a finitely generated ideal, and assume that the set
    \[
    Z(I) := \{ x \in \Omega \mid \forall f \in I, f(x) = 0 \}
    \]
    is a discrete subset of \( \Omega \). Then \( I \) is a principal ideal; i.e., there exists \( f \in \mathcal{E}^\infty(\Omega) \) such that \( I = (f) \).
\end{theo}
{\bf{Proof.}}
Let \( I = (f_1, f_2, \dots, f_n) \subset \mathcal{E}^\infty(\Omega) \). Assume that the common zero set
\[
Z(I) = \{ x \in \mathbb{R} \mid f_1(x) = \dots = f_n(x) = 0 \}
\]
is discrete. Since \( Z(I) \) is closed and has no accumulation points, we can enumerate its elements as \( Z(I) = \{ x_k \}_{k \in \mathbb{Z}} \) (possibly finite).\\
For each point \( x_k \in Z(I) \), choose an open neighborhood \( U_k \) of \( x_k \) such that:
\begin{itemize}
    \item \( U_k \cap Z(I) = \{ x_k \} \),
    \item \( \overline{U_k} \cap \overline{U_j} = \emptyset \) for \( j \neq k \),
    \item On each \( U_k \), the restriction \( I|_{U_k} \subset \mathcal{E}^\infty(U_k) \) is principal.
\end{itemize}
The last point holds because in dimension one, finitely generated smooth function ideals are always locally principal when their common zeros are isolated.
Thus, there exists a smooth function \( g_k \in \mathcal{E}^\infty(U_k) \) such that
\[
I|_{U_k} = (g_k) \subset \mathcal{E}^\infty(U_k).
\]

Let \( \phi_k \in \mathcal{E}^\infty(U_k) \) be a smooth bump function such that \( \phi_k \equiv 1 \) in a neighborhood of \( x_k \) and $supp\,\phi_k \subset U_k$ is compact. Define:
\[
h_k := \phi_k g_k \in \mathcal{E}^\infty(\Omega),
\]
 Since \( g_k \in I|_{U_k} \), there exist \( h_{k_i} \in C^\infty(U_k) \) such that
\[
g_k = \sum_{i=1}^n h_{k_i} f_i|_{U_k}.
\]
Multiplying both sides by \( \phi_k \) gives:
\[
h_k = \sum_{i=1}^n (\phi_k h_{ki}) f_i,
\]
and since \( \phi_k h_{ki} \in C^\infty(\mathbb{R}) \), we conclude that \( h_k \in I \).

Now, outside of the zero set \( Z(I) \), for any \( x \notin Z(I) \), there exists at least one \( f_i \) such that \( f_i(x) \neq 0 \). So in a neighborhood of \( x \), \( f_i \) generates \( I \), and again \( I \) is locally principal.

Thus, we cover \( \mathbb{R} \) by a locally finite open cover \( \{ V_\alpha \} \), each with a smooth local generator \( g_\alpha \in I \). We then use a smooth partition of unity \( \{ \psi_\alpha \} \) subordinate to this cover and define the global function:
\[
f := \sum_\alpha \psi_\alpha g_\alpha.
\]
This sum is locally finite and therefore defines a smooth function on \( \mathbb{R} \). Since each \( g_\alpha \in I \) and \( \psi_\alpha \in \mathcal{E}^\infty(\mathbb{R}) \), it follows that \( f \in I \). Therefore, \( (f) \subset I \).

To show the reverse inclusion, consider any \( h \in I \). Locally, \( h = \eta_\alpha f \) for some smooth function \( \eta_\alpha \). The local expressions can be patched globally using a partition of unity to define a global smooth function \( \eta \in \mathcal{E}^\infty(\mathbb{R}) \) such that \( h = \eta f \). Thus \( h \in (f) \), and so \( I \subset (f) \).

In dimension one, we give a converse of Corollary 5.
\begin{cor}
    Let $I\subset\mathcal{E}^\infty (\Omega)$ be a finitely generated $L$ojasiewicz ideal. Then $I$ is closed in the Frechet topolofy of $\mathcal{E}^\infty (\Omega)$.
\end{cor}
{\bf{Proof.}}
By Corollary 7, the set of zeros of the ideal $I$ is discrete, hence the ideal $I$ is principal, by Theorem 10. Then there exists $f\in I$ such that $ I = (f)\mathcal{E}^\infty (\Omega )$. Let $g\in\overline{I}$, then $\frac{h}{f}\in \mathcal{E}^\infty (\Omega )$ so there exists $g\in \mathcal{E}^\infty (\Omega )$ with $ h= gf$, hence $h\in I$.
\bibliographystyle{alpha}
\bibliography{sample}
{\textbf{Abdelhafed El Khadiri\\
Universty Ibn Tofail\\ Department of Mathematics\\
elkhadiri.abdelhafed@uit.ac.ma}}

\end{document}